\documentclass[reqno]{amsart} 
\usepackage{amsmath}
\usepackage{paralist}
\usepackage{amssymb}
\usepackage{breqn}

\usepackage{titlesec}
\usepackage[misc]{ifsym}
\usepackage{epsfig} 
\usepackage{epstopdf} 
\usepackage[colorlinks=true]{hyperref}
   
\hypersetup{urlcolor=blue, citecolor=red}

\allowdisplaybreaks

\titleformat{\section}
  {\normalfont\bfseries\centering}
  {\thesection.}{1em}{}

\usepackage[a4paper, left=1in, right=1in, top=1in, bottom=1in]{geometry}

\theoremstyle{definition}

\newcommand{\ie}{i.e.}

\title[Singular Superlinear Eqns on Ext Domains]{Existence and nonexistence of sign-changing solutions for linearly perturbed superlinear equations on exterior domains} 

\author [ Md Suzan Ahamed, Joseph Iaia ]{}



\begin{document}

\maketitle

\centerline{\scshape
Md Suzan Ahamed$^{{\href{mailto:mdsuzanahamed@my.unt.edu}{\textrm{\Letter}}}*1},$ \scshape Joseph Iaia$^{{\href{mailto:iaia@unt.edu}{\textrm{\Letter}}}*1}$}

\medskip

{\footnotesize \centerline{$^1$Department of Mathematics, University of North Texas, P.O. Box 311430, Denton, 76203, Texas, USA} }

\begin{abstract} 
\noindent In this paper, we study radial solutions of $\Delta u + K(|x|) f(u)+\frac{ (N-2)^2 u}{|x|^{2+(N-2)\delta}} =0, \   0<\delta<2$ in the exterior of the ball of radius $R>0$ in ${\mathbb R}^{N}$ where $f$ grows superlinearly at infinity and is singular at $0$ with $f(u) \sim -\frac{1}{|u|^{q-1}u}$ and $0<q<1$ for small $u$. We assume $K(|x|) \sim |x|^{-\alpha}$ for large $|x|$ and establish the existence of an infinite number of sign-changing solutions when $N+q(N-2) <\alpha <2(N-1).$   We also prove nonexistence for $0<\alpha \leq2$.
 
\end{abstract}

\def\e{\epsilon}

\section{Introduction}
 In this paper, we are interested in radial solutions of
\begin{equation}   \Delta u + K(|x|) f(u)+\frac{(N-2)^2 u}{r^{2+(N-2)\delta}}=0 \quad   0<\delta<2  \textrm{ on } {\mathbb R}^N\backslash B_{R}\label{DE22}  \end{equation}
\begin{equation}  u = 0 \textrm{ on } \partial B_R, \  u\to 0 \textrm{ as } |x| \to \infty \end{equation}
when $N>2$ and where $B_R$ is the ball of radius $R>0$ centered at the origin.

Assuming $u(x) =u(|x|) = u(r)$, then this becomes
\begin{equation}  u'' + \frac{N-1}{r} u' + K(r) f(u)+\frac{ (N-2)^2 u}{r^{2+(N-2)\delta}} =0  \quad  \textrm{for } R<r < \infty, \label{DE1}
\end{equation}
\begin{equation}   u(R) = 0, \lim_{r \to \infty} u(r) = 0.     \label{DE100}   \end{equation}

Numerous papers have proved the existence of {\it positive} solutions of these equations with various nonlinearities $f(u)$ and for various functions $K(|x|) \sim |x|^{-\alpha}$ with  $\alpha>0$. See, for example,  \cite{A},  \cite{BL}-\cite{BL2}, \cite{LSS}-\cite{S}. 
\vskip .1 in
Here, we prove the existence of multiple solutions including {\it sign-changing} solutions for this equation.  We have also proved the existence of {\it sign-changing} solutions in other recent papers: \cite{AI}-\cite{AI2}, \cite{I6}-\cite{I8}.

\vskip .1 in

We assume $f: {\mathbb R} \backslash\{0\} \to {\mathbb R}$ is odd, locally Lipschitz, and
\begin{equation} \tag{\it{H1}}     f(u) = |u|^{p-1}u + g(u)   \textrm{ with } p>1 \textrm{  for large }  |u|
\textrm{ and}  \lim_{u \to \infty} \frac{ |g(u)| }{|u|^{p}} =0.
\end{equation}
We also assume there exists a locally Lipschitz  $g_1: {\mathbb R} \to {\mathbb R}$  such that
\begin{equation} \tag{\it{H2}}     f(u) = -\frac{1}{|u|^{q-1}u} + g_1(u)   \textrm{ with } 0<q<1 \textrm{  for small }  |u| \textrm{ and }\  g_1(0)=0.
\end{equation}
In addition, we assume
\begin{equation} f \textrm{  has a unique positive zero, } \beta, \textrm{ with } f<0 \textrm{  on }(0, \beta) \textrm{ and }  f>0 \textrm{  on } (\beta, \infty).\tag{{\it H3}}  \end{equation} 
Let $F(u) = \int_{0}^{u} f(t) \, dt.$ Since $f$ is odd,  $F$ is even. Also, since $0<q<1$ by{\it (H2)}, it follows that $F$ is continuous with $F(0)=0$. Also,
\begin{align} F \textrm{ has a unique positive zero, }  \gamma, \textrm{ where } 0< \beta< \gamma  \textrm{ with } F<0 \textrm{ on } (0, \gamma) \ \textrm{ and }  F>0 \textrm{  on } (\gamma, \infty).  \label{boss} 
\end{align}
In addition, it follows from {\it (H3)} that there exists $F_0>0$ and $f_0>0$ such that
\begin{equation} F(u)\geq -F_0    \textrm{ and } |u|^{q-1}uf(u) \geq -f_0 \textrm{ for all } u \in {\mathbb R\backslash \{0\}}. \label{boss2} \end{equation}
We also assume $K>0$ and $K'$ are continuous on $[R, \infty).$ In addition, we assume
\begin{align}
\frac{rK'}{K} > -2(N-1),  \ \text{there exist positive } K_0 \ \text{and} \ K_1  \ \text{such that} \ \frac{K_0}{r^{\alpha}} \leq  K\leq \frac{K_1}{r^{\alpha_1}} \notag \\    \textrm{ with }  \ N + q(N-2)<\alpha_1 \leq \alpha<2(N-1)  \textrm{ on } [R, \infty). \tag{\it{H4}} 
\end{align}

{\bf REMARK 1:}
The reason for our interest in this problem with the extra term $\frac{(N-2)^2 u}{r^{2+(N-2)\delta}}$ is that we will see in section 2 that after a change of variables we will ultimately be examining the following boundary problem
\begin{equation}  v'' + h(t) f(v) + \frac{v}{t^{2-\delta}} =0 
 \label{DE10} \end{equation}
with
\begin{equation} v(0)=0=v(R^{2-N}).  \label{DE11} \end{equation}
Now if we examine the linear part of this problem, we have
\begin{equation} v'' + \frac{v}{t^{2-\delta}} =0.  \label{DE12} \end{equation} 

When $\delta>0 $
it is well-known that all solutions of (\ref{DE12}) are continuous and have only a {\it finite number of zeros} on $[0, R^{2-N}]$. In fact, the general solution of (\ref{DE12}) is
\begin{equation} y = c_1 t^{\frac{1}{2}} 
J_{\frac{1}{\delta}} \left( \frac{2}{\delta} t^{\frac{\delta}{2}} \right) + c_2 t^{\frac{1}{2}} Y_{\frac{1}{\delta}}
\left( \frac{2}{\delta} t^{\frac{\delta}{2}}\right)  \label{DDAA} \end{equation}
where $J_{\frac{1}{\delta}}$  and $Y_{\frac{1}{\delta}}$ are the Bessel functions of order $\frac{1}{\delta}$ of the first and second kind, and $c_1$, $c_2$ are constants. 

\vskip .1 in

However, when $\delta=0$ the solution of
\begin{equation}
    v'' + \frac{v}{t^2} =0 \label{linear0001}
\end{equation} is
$$ \ \ \ \ v_a =      c_1  t^{\frac{1}{2}}\cos\left( \frac{\sqrt{3}}{2}\ln t\right)  + c_2 t^{\frac{1}{2}} \sin\left( \frac{\sqrt{3}}{2}\ln t\right). $$

Thus we see that the solutions are continuous and have  {\it infinitely many zeros} on $[0, R^{2-N}]$.
\vskip .1 in

So then ultimately the challenge for us is what happens to solutions of $(\ref{DE10})$-$(\ref{DE11})$ in the case when $ 0<\delta < 2$? Are there solutions of $(\ref{DE10})$-$(\ref{DE11})$ that have infinitely many zeros on $(0, R^{2-N})$ as there are for $(\ref{linear0001})$  or are there infinitely many solutions of $(\ref{DE10})$-$(\ref{DE11})$ where each solution has only a finite number of zeros on $(0, R^{2-N})$? We were able to show it is this second situation which occurs.

 \vskip .1 in
 In this paper, we prove the following: 
\vskip .1 in

{\bf  THEOREM 1:} {\it  Let $N>2$ and assume (H1)-(H4). Then there exists $R_0>0$ and a non-negative integer $n_0$ so that if $R\geq R_0$ then there are infinitely many $C^1[R, \infty)$ solutions, $u_n,$ of (\ref{DE1})-(\ref{DE100}) such that $u_n$ has exactly $n$ zeros on $(R,\infty)$ for each non-negative integer $n\geq n_0.$} 

\vskip .1 in
A similar problem but without the additional linear term was investigated in earlier papers. 
In \cite{AI} we proved Theorem 1 in the case  $  \alpha  > 2(N-1).$ 
In \cite{I8} we proved Theorem 1 in the case $ N+ q(N-2) < \alpha < 2(N-1) $ and we proved nonexistence for $0 < \alpha < N+q(N-2).$

\vskip .1 in
The addition of the linear term adds a number of difficulties to the proofs of several of the lemmas in Section 2. A number of integrals that arise are only defined when $0<\delta <2$.  In addition, the constant $\frac{1}{\delta}$ appears as an upper bound in a number of estimates.

\vskip .2 in

{\bf REMARK 2:}  Solutions of (\ref{DE1})-(\ref{DE100}) have continuous second derivatives except at points where $u(r_0)=0$  due to the fact that  $\lim\limits_{u \to 0} |f(u)| =\infty$. Solutions, however, do turn out to be $C^{1}[R, \infty)$. Therefore, by a {\it $C^1[R, \infty)$ solution of (\ref{DE1})-(\ref{DE100})} we mean  $u\in C^{1}[R, \infty)$ such that $\frac{1}{|u|^q}$ is integrable on $[R, \infty)$ and $$r^{N-1} u' + \int_{R}^r t^{N-1} K f(u)+\int_R^r \frac{u(t)}{t^{3-N+(N-2)\delta}}\ dt = R^{N-1}u'(R) \quad \text{for} \ r\geq R $$
$$ u(R) = 0, \ \ \lim_{r \to \infty} u(r) = 0.$$

\vskip .1 in
In this paper, we also consider the case where
\begin{equation} 
-\frac{rK'}{K} <\alpha    \textrm{ with } 0 <\alpha \leq 2  \textrm{ on } [R, \infty). \tag{\it H5}  
\end{equation}

We also prove:
\vskip .1 in

{\bf  THEOREM 2:} {\it Assume {\it (H1)-(H3)} and {\it (H5)}. If $R>0$ and $N>2,$  then there are no   $C^1[R, \infty)$ solutions  of (\ref{DE1})-(\ref{DE100}) when $0<\alpha \leq 2$. }   
\vskip.2 in

\section {Preliminaries for Theorem 1}

Let $u$ solve $(\ref{DE1})$-$(\ref{DE100})$ and suppose {\it(H1)}-{\it (H4)} are satisfied. Now let $u(r) = v(r^{2-N}).$ Then we see (\ref{DE1})-(\ref{DE100}) is equivalent to

\begin{equation}
v''(t) + h(t) f(v(t))+\frac{ v(t)}{t^{2-\delta}} = 0  \   \textrm{ for } 0<t< R^{2-N},  \label{e9} \end{equation}
\begin{equation} v(0) =0, \  v(R^{2-N})=0 \label{e10}  \end{equation}

where \begin{equation}h(t) = \frac{t^{\frac{2(N-1)}{2-N}}K(t ^{\frac{1}{2-N}})}{(N-2)^2}>0,  \ 0<\delta<2 . \label{h eq}\end{equation}
It follows from {\it (H4)} that $h'<0.$  

\vskip .1 in
We now define $\tilde \alpha = \frac{2(N-1)-\alpha}{N-2}$ and $\tilde \alpha_1 = \frac{2(N-1)-\alpha_1}{N-2}.$ Note that the assumption $N + q(N-2)<\alpha_1 \leq \alpha<2(N-1)$ from {\it (H4)} implies
\begin{equation}
    0<\tilde\alpha + q \leq  \tilde \alpha_1 +q <1. \label{alphatilde}
\end{equation}
It also follows from {\it (H4)} and  (\ref{h eq}) that there exist \ $h_0>0 \ \text{and } \  h_1>0$ such that
\begin{equation}    \  h_0 t^{-\tilde \alpha} \leq h \leq h_1 t^{-\tilde \alpha_1}  \textrm{ for }  0< t< R^{2-N},  \label{joni} 
\end{equation}
and from this and $(\ref{alphatilde})$ it follows that $\frac{h(t) }{t^q}$ and $h(t)$  are integrable on $(0, R^{2-N}].$
\vskip .1 in
We now consider the following initial value problem
\begin{equation} 
v_a''(t) + h(t) f(v_a(t))+\frac{ v_a(t)}{t^{2-\delta}} = 0   \textrm{ for } t>0, \label{e7} 
\end{equation}
\begin{equation} v_a(0) =0, \ v_a'(0)=a> 0. \label{e8}  \end{equation}
{\bf REMARK 3:} We do not consider the case when $\delta=0$ with the initial conditions in $(\ref{e8})$ because if there is such a solution then we see that $tv_a''\sim -a$ for small $t$ and so then dividing by $t$ and integrating gives that $v_a'\sim -a \ln(t)$ which blows up as $t$ goes to $0$ and so we cannot have $v_a(0)=0$ and $v_a'(0)=-a$ when $\delta =0$.

\vskip .1 in
{\bf Lemma 2.1:} Let  $N>2$,  $a>0,$  and assume {\it(H1)}-{\it (H4)}. Then there exists an $\e >0$ such that there is a unique solution of  (\ref{e7})-(\ref{e8}) on $[0, \e].$

\vskip .1 in

{\bf Proof:} We first assume $v_a$ solves $(\ref{e7})$-$(\ref{e8})$, then integrate (\ref{e7}) on $(0, t)$ and use (\ref{e8}) along with {\it (H2)} to obtain
\begin{equation} v_a' = a - \int_0^{t} h(x)f(v_a(x))  \, dx - \int_0^{t} x^{\delta-2} v_a(x) \ dx . \label{thomas} \end{equation}
Now, integrate (\ref{thomas})  on $(0,t)$ to obtain
\begin{equation} v_a = at -\int_0^{t}\int_0^{s} h(x)f(v_a(x))  \, dx \, ds - \int_0^t\int_0^{s} x^{\delta-2} v_a(x) \ dx \ ds .  \label{jeter}  \end{equation}
\vskip .1 in
Next, let $w(t)=\frac{v_a(t)}{t}.$ Then, we obtain $w(0)=a>0$ and
\begin{equation} 
w(t) = a - \frac{1}{t}\int_0^{t}\int_0^{s} h(x)f(xw(x))  \, dx \, ds - \frac{1}{t} \int_0^t\int_0^{s} x^{\delta-1} w(x) \ dx \ ds \quad \textrm{ for } t>0.  \label{w eqn}  \end{equation}
Now using the contraction mapping principle, we prove the existence of a solution of  (\ref{e7})-(\ref{e8}) on $[0, \epsilon]$ for some $\epsilon>0.$ Let
$$ \mathcal{B}:=\{ w\in C[0, \epsilon] \ | \ w(0)=a \ \text{and} \  \left|w(t)-a\right|\leq \frac{a}{2} \ \text{on} \  [0, \epsilon]\},$$ where $C[0,\epsilon]$  is the set of all continuous functions on $[0,\epsilon],$ $\epsilon>0,$ and
$$||w||=\sup_{x\in [0,\epsilon]} |w(x)|.$$ Then $(\mathcal{B},||.||)$ is a Banach space \cite{book2010}. Now define a map $T$ on $\mathcal{B}$ by: 
$$Tw(t)= \begin{cases}
      a    \ \ \quad\quad\quad \quad \quad \quad \quad \ \ \ \ \ \ \  \ \ \ \ \ \ \ \  \ \ \ \ \ \ \  \  \quad  \quad\quad\quad  \quad\quad\quad\quad\quad    \quad   \ \   \          \textrm { for } t=0\\
      a - \frac{1}{t}\int_0^{t}\int_0^{s} h(x)f(xw(x))  \, dx \, ds - \frac{1}{t} \int_0^t\int_0^{s} x^{\delta-1} w(x) \ dx \ ds  \quad \quad \textrm{\ for } \   0<t \leq \epsilon. 
         \end{cases}$$
Next, note that for any $w\in \mathcal{B},$ we have
\begin{equation}
    0<\frac{a}{2}\leq w(x)\leq \frac{3a}{2} \quad \text{on} \ [0, \epsilon].\label{wb001}
\end{equation}
Then $\left|\frac{-1}{x^qw^q(x)}\right|\leq \frac{2^q x^{-q}}{a^q} \ \text{on} \  (0, \epsilon ].$ Since $g_1$ is locally Lipschitz by $(\it{H2})$ and $g_1(0)=0,$ therefore there exists $L_1>0$ such that for some $\epsilon_1>0$
\begin{equation}
    \left|g_1(xw(x))\right|\leq \frac{3aL_1}{2} x \quad \text{on} \  \left [ 0, \epsilon_1 \right ]. \label{g1}
\end{equation}
Let $\epsilon_2:=\min \{\epsilon, \epsilon_1\}.$ Then it follows from this, $(\ref{joni})$, and $(\ref{g1})$
$$|h(x)f(xw(x))|\leq \frac{2^qh_1}{a^q}x^{-(\tilde{\alpha}+q)}+\frac{3ah_1L_1}{2}x^{1-\tilde{\alpha}} \quad \text{on} \  \left ( 0, \epsilon_2 \right ].$$
Integrating on $(0, t)\subset (0, \epsilon_2)$, and using $(\ref{alphatilde})$,  we obtain
\begin{equation}
    \int_0^t |h(x)f(xw(x))| \ dx \leq C_1 t^{1-\tilde{\alpha}-q}+C_2 t^{2-\tilde{\alpha}}\to 0 \quad \text{as} \ t\to 0^+, \label{hf1}
\end{equation}
where $C_1:=\frac{2^qh_1}{a^q(1-\tilde{\alpha}-q)}$ and $C_2:=\frac{3ah_1L_1}{2(1-\tilde{\alpha})}.$ \\
Integrating $(\ref{hf1})$ again, we obtain
\begin{equation}
    \frac{1}{t}\int_0^{t}\int_0^{s} |h(x)f(xw(x))|  \, dx \, ds \to 0 \quad \text{as} \ t\to 0^+.\label{hf3}
\end{equation}
Since $|x^{\delta-1}w(x)|\leq \frac{3a}{2}x^{\delta-1}$ on $(0, \epsilon_2)$ $\ (\because w\in \mathcal{B})$, integrating both sides of this inequality on $(0, t)\subset (0, \epsilon_2)$ gives
$$\int_0^t |x^{\delta-1}w(x)|\ dx \leq \frac{3a}{2\delta}t^\delta.$$
Integrating again on $(0, t)\subset (0, \epsilon_2)$, we obtain
\begin{equation}
\frac{1}{t}\int_0^t \int_0^s |x^{\delta-1}w(x)|\ dx \ ds \leq \frac{3a}{2\delta(\delta+1)}t^\delta \to 0 \quad \text{as} \ t\to 0^+.  \label{hf4}
\end{equation}
We now show that $Tw\in \mathcal{B}$ for each $w \in \mathcal{B}$ if $\epsilon_2>0$ is sufficiently  small. For any $w\in \mathcal{B}$, $Tw(0)=a$ and it is clear from the definition and $(\ref{hf3})$-$(\ref{hf4})$ that $Tw(t)$ is continuous for any $t\in (0, \epsilon].$ It is also clear from $(\ref{hf3})$ that
$$ |Tw(t)-a|\leq \frac{1}{t}\int_0^{t}\int_0^{s}|h(x)f(xw(x))|  \, dx \, ds + \frac{1}{t}\int_0^t \int_0^s |x^{\delta-1}w(x)|\ dx \ ds \ \to 0 \quad \text{as} \ t\to 0^+.$$
Thus $|Tw(t)-a|\leq \frac{a}{2}$ for sufficiently small $t>0.$ Therefore $T:\mathcal{B}\longrightarrow \mathcal{B}\text{ if } \epsilon$ is sufficiently small.
\vskip .1 in
We next prove that $T$  is a contraction mapping if $\epsilon$ is sufficiently small. Let $w_1,w_2 \in \mathcal{B}$. Then
\begin{equation}    |Tw_1- Tw_2|\leq \frac{1}{t}\int_0^{t}\int_0^{s} |h(x)[f(xw_1)-f(xw_2)]|  \, dx \, ds + \frac{1}{t}\int_0^t \int_0^s |x^{\delta-1}(w_1-w_2)|\ dx \ ds. \label{tw1}
\end{equation}
Note that by {\it (H2)}
$$| f(xw_1)-f(xw_2)|\leq \frac{1}{x^q} \left|\frac{1}{{w_1}^q}-\frac{1}{{w_2}^q}\right|+L_1  x |w_1-w_2|,$$
where $L_1$ is the Lipschitz constant for $g_1.$ Then using the mean value theorem for any $x\in (0, \epsilon)$ there is a $c_x$ with $\frac{a}{2}\leq c_x \leq \frac{3a}{2}$ such that
$$ \frac{1}{x^q}\left|\frac{1}{{w_1}^q}-\frac{1}{{w_2}^q}\right|=\frac{q}{x^q c_x^{q+1}}|w_1-w_2|\leq \frac{2^{q+1}q}{a^{q+1}} |w_1-w_2|x^{-q}.$$
It follows from this, $(\ref{alphatilde}), (\ref{joni})$, and integrating on $(0, t)$ that
\[\int_0^t|h(x)(f(xw_1)-f(xw_2))| \ dx\leq \frac{2^{q+1}q h_2}{(1-\tilde{\alpha}-q)a^{q+1}} |w_1-w_2|t^{1-(\tilde{\alpha}+q)}+\frac{L_1 h_2}{2-\tilde{\alpha}} |w_1-w_2|t^{2-\tilde{\alpha}}.\]
Integrating again on $(0, t)$, we obtain
\[\int_0^t \int_0^s|h(x)(f(xw_1)-f(xw_2))| \ dx \ ds\leq C_3 |w_1-w_2|t^{2-(\tilde{\alpha}+q)}+ C_4|w_1-w_2|t^{3-\tilde{\alpha}},\]
where $C_3=\frac{2^{q+1}q h_1}{((1-\tilde{\alpha}-q))(2-\tilde{\alpha}-q)a^{q+1}}$ and $C_4=\frac{L_1 h_1}{(2-\tilde{\alpha})(3-\tilde{\alpha})}.$ Thus
\begin{equation}
    \frac{1}{t}\int_0^t \int_0^s|h(x)(f(xw_1)-f(xw_2))| \ dx \ ds\leq (C_3 \epsilon^{1-(\tilde{\alpha}+q)}+ C_4 \epsilon^{2-\tilde{\alpha}})||w_1-w_2|| \quad \text{on} \ (0, \epsilon). \label{nr1}
\end{equation}
Also, since $|x^{\delta-1}(w_1-w_2)|\leq x^{\delta-1}||w_1-w_2||$, then integrating twice on $(0, t)\subset (0, \epsilon),$ we obtain
\begin{equation}
    \frac{1}{t} \int_0^t \int_0^s |x^{\delta-1}(w_1-w_2)| \ dx \ ds \leq \frac{\epsilon^\delta}{\delta (\delta+1)}  ||w_1-w_2||. \label{xw1}
\end{equation}
Using $(\ref{nr1})$ and $(\ref{xw1})$ in $(\ref{tw1})$, we see that
$$|Tw_1- Tw_2|\leq C_\epsilon ||w_1-w_2||,$$
where $C_\epsilon :=C_3 \epsilon^{1-(\tilde{\alpha}+q)}+ C_4 \epsilon^{2-\tilde{\alpha}}+\frac{\epsilon^\delta}{\delta (\delta+1)}.$ Note that $C_\epsilon\to 0$ as $\epsilon\to 0^+.$ Therefore, $0<C_\epsilon<1$ for sufficiently small $\epsilon>0$, and hence $T$ is a contraction on $\mathcal{B}.$ Then by the contraction mapping principle there exists a unique $w\in \mathcal{B}$ such that $Tw=w$ on $[0, \epsilon]$ for some $\epsilon>0$, and hence $v_a(t)=tw(t)$ is a solution of $(\ref{e7})$-$(\ref{e8})$ on $[0, \epsilon].$ \\
Moreover, from $(\ref{wb001}),$ we see that $\frac{a}{2}t\leq v_a\leq \frac{3a}{2}t.$ So $v_a$ is bounded on $[0, \epsilon].$ Also, from $(\ref{hf3})$ and $(\ref{hf4})$ it follows that $v_a'$ is bounded on $[0, \epsilon].$ This completes Lemma 2.1. \qed
\vskip .1 in

{\bf Lemma 2.2:}  Let $N>2$ and assume {\it (H1)}-{\it (H4)}. Suppose $v_a$ solves (\ref{e7})-(\ref{e8}) with $a>0$ on some maximal half-open interval $[0,d)$ with $\epsilon\leq d \leq R^{2-N}$. Then there exists a constant $C>0$ where $C$ is independent of $d$  such that $|v_a(t)|\leq C$ and $|v_a'(t)| \leq C$ on $[0, d]\subset [0, R^{2-N}]$. In particular, it follows from this that $v_a$ and $v_a'$ are continuous, can be extended to all of $[0, R^{2-N}]$, and that $v_a$ varies continuously with $a> 0$.

\vskip .1 in

{\bf Proof:}   Note it follows from $(\ref{e7})$-$(\ref{e8})$ that
$$\left(\frac{1}{2}v_a'{^2}(t)+\frac{ v_a^2(t)}{2t^{2-\delta}}+h(t)F(v_a(t))\right)'\leq h'(t)F(v_a(t)).$$
Integrating both sides on $[\epsilon, t]$ where $[0, \epsilon]$ is the interval of existence obtained in Lemma 2.1, we obtain
$$\frac{1}{2}v_a'{^2}(t)+h(t)F(v_a(t))\leq \frac{1}{2}v_a'{^2}(\epsilon)+\frac{ v_a^2(\epsilon)}{2\epsilon^{2-\delta}}+h(\epsilon)F(v_a(\epsilon))+\int_{\epsilon}^t h'(s)F(v_a(s)) \ ds.$$
It follows from $(\ref{boss2})$ that $F\geq-F_0$ for some $F_0.$ Since $h>0$, therefore $hF\geq-hF_0, \ h'F\leq-F_0h' \ (\because h'<0 \ \text{from the comment after (\ref{h eq})})$, and hence
$$v_a'{^2}(t)\leq v_a'{^2}(\epsilon)+\frac{ v_a^2(\epsilon)}{\epsilon^{2-\delta}}+2h(\epsilon)(F(v_a(\epsilon))+F_0) \quad \text{on} \ [\epsilon, d].$$
Taking the square root on both sides, we obtain $|v_a'(t)|\leq C_5 \ \text{on} \ [\epsilon, d]$, where
$$ C_5:=\sqrt{v_a'{^2}(\epsilon)+\frac{ v_a^2(\epsilon)}{\epsilon^{2-\delta}}+2h(\epsilon)(F(v_a(\epsilon))+F_0)}.$$
Also, since $v_a(t)=v_a(\epsilon)+\int_{\epsilon}^t v_a'(s) \ ds$, then $|v_a(t)|\leq C_6 \ \text{on} \ [\epsilon, d]$ where $C_6:=C_5(R^{2-N}-\epsilon)+|v_a(\epsilon)|.$
\vskip .1 in
Now combining this with the bounds obtained for $v_a$ and $v_a'$ on $[0, \epsilon]$ in Lemma 2.1, we see that $v_a$ and $v_a'$ are uniformly bounded on $[0,d]$ by a constant that only depends on $\epsilon$ and $R.$ It follows from this that $v_a$ and $v_a'$ can be extended to be  defined on $[0, R^{2-N}]$,  and $v_a$, $v_a'$ are uniformly bounded on $[0, R^{2-N}].$ It also follows from this that $v_a$ varies continuously with $a>0$. This completes the proof. \qed

\vskip .1 in

{\bf Lemma 2.3:} Let $N>2$,  and assume {\it (H1)}-{\it (H4)} or {\it (H1)}-{\it (H3)}, and {\it (H5)}. Let $v_a$ solve (\ref{e7})-(\ref{e8}). If $v_a(z_a)=0$ and $M_a$ is a local maximum (or local minimum) with $0\leq z_a< M_a$ and $v_a'\geq 0$ (or $v_a'\leq 0$) on $(z_a, M_a)$, then there is a $C_7>0$ such that $|v_a(M_a)| \geq C_7$ where $C_7$ is independent of $a.$ In addition, there exists $R_0>0$ such that $|v_a(M_a)| >\gamma$  if $R>R_0.$

 \vskip .1 in
 
{\bf Proof:} Let us suppose $v_a>0$ and $v_a'\geq 0 $ on $ (z_a, M_a)$ . If $v_a(M_a) >\gamma$ then we are done and we can just use $C_7 =\gamma$. So let us suppose $ 0< v_a(M_a) \leq \gamma$.
Then rewrite $(\ref{e7})$ as
$$ ((t-z_a)v_a' - v_a)' + (t-z_a) h g(v_a) + (t-z_a) \frac{v_a}{t^{2-\delta}} = \frac{(t-z_a)h}{v_a^q}. $$
Integrate on $(z_a,t)$ to obtain
$$ (t-z_a) v'_a-v_a + \int_{z_a}^{t} (s-z_a)h g(v_a) \, ds + \int_{z_a}^{t} (s-z_a)\frac{v_a}{s^{2-\delta}} \, ds= \int_{z_a}^{t} \frac{(s-z_a)h}{v_a^q} \,ds . $$
Now let $t=M_a$ to get
\begin{equation}
     \int_{z_a}^{M_a} (s-z_a)h g(v_a) \ ds + \int_{z_a}^{M_a} (s-z_a)\frac{v_a}{s^{2-\delta}}\, ds = \int_{z_a}^{M_a} \frac{(s-z_a)h}{v_a^q} \, ds + v_a(M_a). \label{456} 
\end{equation}
Since $v_a$ is increasing on $(z_a, M_a),$ then we see
$$ \int_{z_a}^{M_a} (s-z_a)\frac{v_a}{s^{2-\delta}} \, ds \leq \int_{z_a}^{M_a} \frac{v_a}{s^{1-\delta}} \, ds \leq \frac{v_a(M_a)(M_a^\delta - z_a^\delta)}{\delta} $$
and
\begin{equation}
    \int_{z_a}^{M_a} \frac{(s-z_a)h}{v_a^q} \, ds \geq \frac{1}{v_a^q(M_a)} \int_{z_a}^{M_a} (s-z_a)h \, ds. \label{aka001}
\end{equation}
Since $v_a(M_a) \leq \gamma,$ then estimating the leftmost term of $(\ref{456})$ we obtain
\begin{equation}
    \int_{z_a}^{M_a} (s-z_a)h g(v_a) \, ds \leq \left( \max_{[0, \gamma]} |g| \right) \int_{z_a}^{M_a} (s-z_a)h \, ds.\label{aka002}
\end{equation}
Combining $(\ref{aka001})$-$(\ref{aka002})$ and inserting into $(\ref{456})$ gives
\begin{equation}
\left( \max_{[0, \gamma]} |g| \right) \int_{z_a}^{M_a} (s-z_a)h \, ds \geq \frac{1}{v_a^q(M_a)} \int_{z_a}^{M_a} (s-z_a)h \, ds + v_a(M_a)\left( 1 - \frac{(M_a^\delta - z_a^\delta)}{\delta} \right) \label{457}.
\end{equation}
{\bf Case 1:} $\left(1 - \frac{(M_a^\delta - z_a^\delta)}{\delta}\right) \geq 0$\\
In this case, it follows from $(\ref{457})$ that
\begin{equation}
    v_a(M_a) \geq \frac{1}{ \left(\max\limits_{[0, \gamma]} |g|\right)^{\frac{1}{q}} }. \label{aka003}
\end{equation}
{\bf Case 2:} $\left(1 - \frac{(M_a^\delta - z_a^\delta)}{\delta}\right) < 0$\\
 In this case $(\ref{457})$ becomes
 \begin{equation}
     \left( \max_{[0, \gamma]} |g|\right ) \int_{z_a}^{M_a} (s-z_a)h \, ds + v_a(M_a)\left( \frac{(M_a^\delta - z_a^\delta)}{\delta} -1 \right)\geq \frac{1}{v_a^q(M_a)} \int_{z_a}^{M_a} (s-z_a)h \, ds \label{458}.
 \end{equation}
Thus,
\begin{equation}
  \left( \max_{[0, \gamma]} |g|\right ) \int_{z_a}^{M_a} (s-z_a)h \, ds + \frac{  R^{(2-N)\delta}\gamma }{\delta}\geq \frac{1}{v_a^q(M_a)} \int_{z_a}^{M_a} (s-z_a)h \, ds. \label{459}  
\end{equation}
Hence
\begin{equation}
   v_a^q(M_a) \geq \frac{1}{\max\limits _{[0, \gamma]}| g| + \frac{ R^{(2-N)\delta}\gamma }{\delta \int_{z_a}^{M_a} (s-z_a)h \, ds }}. \label{460} 
\end{equation}
Finally, we will find a lower bound for $\int_{z_a}^{M_a} (s-z_a)h \, ds .$ Whether we assume {\it (H4)} or {\it (H5)} it follows that $2(N-1)K+rK'\geq 0$ which implies that $h$ is decreasing, and hence we see $$\int_{z_a}^{M_a} (s-z_a)h \, ds \geq h(R^{2-N})\int_{z_a}^{M_a} (s-z_a) \, ds = \frac{h(R^{2-N})(M_a-z_a)^2}{2}. $$
Since we are in case 2, then we see $M_a \geq (z_a^\delta + \delta)^{\frac{1}{\delta}},$ and hence
$$\int_{z_a}^{M_a} (s-z_a)h \, ds \geq \frac{h(R^{2-N})(M_a-z_a)^2}{2}\geq \frac{h(R^{2-N})}{2}\left( (z_a^\delta + \delta)^{\frac{1}{\delta}} - z_a \right)^2. $$
Now note that if we let $b(x) =( x^{\delta} + \delta)^{\frac{1}{\delta}} -x$ where $x \in [0,R^{2-N}]$ and $0< \delta <2,$ then we see $b(x)$ is decreasing for $1<\delta<2$, is constant for $\delta =1$, and is increasing for $0<\delta < 1.$ So in all cases we see $b(x) \geq b_0 = \min\{ \delta^{\frac{1}{\delta}}, (R^{(2-N)\delta} +\delta)^{\frac{1}{\delta}} - R^{2-N})\}>0.$
Thus we see $$\int_{z_a}^{M_a} (s-z_a)h \, ds \geq h(R^{2-N})\int_{z_a}^{M_a} (s-z_a) \, ds = \frac{h(R^{2-N})(M_a-z_a)^2}{2} \geq \frac{b_0^2 h(R^{2-N})}{2}. $$
Then from $(\ref{460})$, we obtain
\begin{equation}
    v_a(M_a) \geq \frac{1}{ \left( \max\limits _{[0, \gamma]}| g| + \frac{2 R^{(2-N)\delta}\gamma}{\delta b_0^2 h(R^{2-N})} \right)^{\frac{1}{q}} }. \label{aka005}
\end{equation}
Finally, combining $(\ref{aka003}), \ (\ref{aka005})$ and the possibility that $v_a(M_a)\geq \gamma,$ we let
$$ C_7 = \min \left \{ \gamma, \frac{1}{ \left( \max\limits _{[0, \gamma]}| g| + \frac{2 R^{(2-N)\delta}\gamma}{\delta b_0^2 h(R^{2-N})} \right)^{\frac{1}{q}} }, \frac{1}{ \left(\max\limits_{[0, \gamma]} |g|\right)^{\frac{1}{q} }} \right \}.$$
Hence $|v_a(M_a)| \geq C_7>0.$ This completes the first part of the lemma.
\vskip .1 in
For the second part of the lemma, we again rewrite $(\ref{e7})$ as
$$ ((t-z_a)v_a' - v_a)' + (t-z_a) h f(v_a) + (t-z_a) \frac{ v_a}{t^{2-\delta}} = 0. $$
Integrating on $(z_a, M_a)$ we obtain
\begin{equation}
    v_a(M_a)=\int_{z_a}^{M_a} (t-z_a)h(t) f(v_a(t)) \ dt + \int_{z_a}^{M_a} (t-z_a)\frac{ v_a(t)}{t^{2-\delta}}\ dt. \label{abc01}
\end{equation}
Note that 
$$\int_{z_a}^{M_a} (t-z_a)\frac{ v_a(t)}{t^{2-\delta}}\ dt\leq \frac{v_a(M_a)}{\delta}\left(M_a^\delta-z_a^\delta\right)\leq \frac{v_a(M_a)}{\delta}R^{(2-N)\delta}.$$
Since $\frac{1}{\delta}R^{(2-N)\delta}\to 0$ as $R\to \infty,$ therefore there exists $R_1>0$ such that 
\begin{equation}
    \int_{z_a}^{M_a} (t-z_a)\frac{ v_a(t)}{t^{2-\delta}}\ dt\leq \frac{ 1}{2} v_a(M_a) \quad \text{for} \ R\geq R_1. \label{abc02}
\end{equation}
Also, using $(\ref{joni})$ and {\it (H1)} this then gives $f(v_a)\leq C_* v_a^p \ \text{for some}\ C_*>0,$ and so we obtain 
\begin{equation}
    \int_{z_a}^{M_a}(t-z_a)h(t) f(v_a(t)) \ dt\leq C_* h_1\int_{z_a}^{M_a} t^{1-\tilde{\alpha_1}} v_a^p(t) \ dt\leq \frac{C_* h_1 v_a^p(M_a)}{2-\tilde{\alpha_1}} R^{(2-N)(2-\tilde{\alpha_1})}. \label{abc03}
\end{equation}
Now using $(\ref{abc02})$ and $(\ref{abc03})$ in $(\ref{abc01})$, we obtain for $R>R_1$
$$ v_a^{p-1}(M_a)\geq \frac{2-\tilde{\alpha_1}}{2h_1}R^{(N-2)(2-\tilde{\alpha_1})}\to \infty \  \text{as} \ R\to \infty.$$
Hence there exists $R_0>R_1$  such that $v_a(M_a)>\gamma \ \text{if} \ R>R_0 \ \text{for all} \ a>0.$ This completes the lemma.\qed

\vskip .1 in

{\bf Lemma 2.4:} Let $N>2, \ a>0$, and assume {\it (H1)}-{\it (H4)}. Suppose $v_a$ solves (\ref{e7})-(\ref{e8}). Then $v_a$ has at most a finite number of zeros on $[0, R^{2-N}].$

\vskip .1 in

{\bf Proof:} Suppose by way of contradiction that $v_a$ has an infinite number of zeros, $z_k$,  and without loss of generality, let us assume $z_k<z_{k+1}.$  Then $v_a$ would have an infinite number of extrema, $M_k$, with $z_k < M_{k}<z_{k+1}$. Since $z_k$ and $M_k$ are in the compact set $[0, R^{2-N}],$ it follows that $z_k$ would have a  limit point $z$, and also we would have $M_k\to z. $  Now, by the mean value theorem and Lemma 2.3,   $C_7\leq |v_a(M_{k})| =|v_a(M_k) - v_a(z_k)| = |v_a'(c_k)| |M_k-z_k|$ for $z_k <c_k < M_k$. However, $|M_k-z_k| \to 0$. Thus $|v_a'(c_k)|\to \infty$, contradicting that the $|v_a'(c_k)|$ are bounded by Lemma 2.2. Thus, $v_a$ has at most a finite number of zeros.\qed

\vskip .1 in

{\bf Lemma 2.5:} Let $N>2,$ $a>0$, and assume {\it (H1)}-{\it (H4)}. Suppose $v_a$ solves (\ref{e7})-(\ref{e8}). Then
$$\max\limits_{[0, \ R^{2-N}]} |v_a| \to \infty \  \text{as} \ a \to \infty.$$
\vskip .1 in

{\bf Proof:} Suppose by way of contradiction that there is a $C_8>0$ such that $|v_a(t)|\leq C_8$ for all $t\in [0, \ R^{2-N}]$ and for all $a>0.$ Since $g_1$ is continuous, there is a $C_{9}>0$ such that $|g_1(v_a(t))|\leq C_{9}$ for all $t\in [0, \ R^{2-N}].$ \\
Now, consider the following initial value problem
\begin{equation} 
y''(t) +\frac{ y(t)}{t^{2-\delta}} = 0   \textrm{ for } t>0, \label{b1} 
\end{equation}
\begin{equation} y(0) =0,\  y'(0)=1. \label{b2}  \end{equation}
Note that this is equivalent to Bessel's equation of order $\frac{1}{\delta}$ and in fact $y(t)=A t^{\frac{1}{2}}J_{\frac{1}{\delta}}(\frac{2}{\delta}t^{\frac{\delta}{2}})$ with $A>0.$
Now multiply $(\ref{e7})$ by $ y(t),$ $(\ref{b1})$ by $ v_a(t)$ and subtract to obtain
\begin{equation}
    (y(t)v_a'(t)-v_a(t)y'(t))'+h(t)y(t)f(v_a(t))=0. \label{vy000001}
\end{equation}
Let $J_{\frac{1}{\delta}, 1}$ be the first positive zero of the Bessel function $J_{\frac{1}{\delta}}$. Using {\it (H2)} we see $h(t)y(t)f(v_a)\leq h(t)y(t)g_1(v_a)$ and integrating $(\ref{vy000001})$ on $(0, t)\subset (0, t_1)$ where $t_1>0$ is sufficiently small and $t_1<J_{\frac{1}{\delta}, 1},$ we obtain
$$y(t)v_a'(t)-v_a(t)y'(t)+\int_0^t h(s)y(s)g_1(v_a(s)) \ ds\geq0.$$
 Then rewriting the above inequality, we get
$$\left( \frac{v_a(t)}{y(t)} \right) '+\frac{1}{y^2(t)} \int_0^t h(s)y(s)g_1(v_a(s)) \ ds\geq 0 \quad \forall \ t\in (0, t_1).$$
Integrating again on $(0, t)\subset (0, t_1),$ we see
\begin{equation}
    v_a(t)+y(t) \int_0^t \frac{1}{y^2(t)} \int_0^s h(x)y(x)g_1(v_a(x)) \ dx \ ds\geq ay(t). \label{vy2}
\end{equation}
Now, since $y(t)\geq0$ on $[0, t_1]$, then from $(\ref{b1})$ we see $y''\leq 0\ \ \text{on } \ [0, t_1].$ Therefore, integrating this twice on $(0, t)$ where $0<t<t_1$, we obtain
\begin{equation}
    y(t)\leq t\ \quad \forall  \ t\in [0, t_1] .\label{y1}
\end{equation}
In addition, since $y(0)=0$ and $y'(0)=1$, there exists a $t_0$ with $0<t_0<t_1$ such that $y(t)\geq \frac{1}{2} t$ for all $t\in [0, t_0].$ Also, since $y(t)>0$ for all $t\in [t_0, t_1)$, then there exists a $C_{y_1}>0$ such that $y(t)\geq C_{y_1}$ on $[t_0, t_1],$ and hence $y(t)\geq C_{y_2}t$ on $[t_0, t_1]$ where $C_{y_2}:=\frac{C_{y_1}}{t_1}.$ Now set $C_y:= \min \{\frac{1}{2}, C_{y_2} \}$. Then $y(t)\geq C_y t$ for all $t\in [0, t_1],$ and using $(\ref{y1})$ we see
\begin{equation}
C_y t \leq y(t)\leq t \quad \forall \ t\in [0, t_1]. \label{y3}
\end{equation}
Using $(\ref{joni})$ and $(\ref{y3})$ in $(\ref{vy2})$, we obtain
$$|ay(t)| \leq \left| v_a(t)+y(t) \int_0^t \frac{1}{y^2(t)} \int_0^s h(x)y(x)g_1(v_a(x)) \ dx \ ds \right|$$
\begin{equation}
    \leq C_8+ \frac{C_{9} h_1 t_1^{2-\tilde{\alpha}}}{C_y^2(1-\tilde{\alpha})(2-\tilde{\alpha})} \quad \text{on} \ [0, t_1]. \label{yc1}
\end{equation}
Now notice that $|ay(t)|\to \infty$ as $a\to \infty$ for any fixed $t$ with $0<t<t_1$. However, the right-hand side of $(\ref{yc1})$ is bounded for all $a$ and the left-hand side is not. This is a contradiction. Thus
\begin{equation}
    \max\limits_{[0, \ t_1]} |v_a| \to \infty \  \text{as} \ a \to \infty. \label{max1}
\end{equation}
Since $\max\limits_{t\in [0, \ R^{2-N}]} |v_a(t)|\geq \max\limits_{t\in [0, \ t_1]} |v_a(t)|,$ therefore we see $\max\limits_{[0, \ R^{2-N}]} |v_a| \to \infty \  \text{as} \ a \to \infty.$\qed

\vskip .1 in

{\bf Lemma 2.6:}   Let $N>2$ and assume {\it (H1)}-{\it(H4)}. Suppose $v_a$ solves (\ref{e7})-(\ref{e8}). Then $v_a$ has a first local maximum, $M_a$,  if $a>0$ is sufficiently large. In addition,  $M_a \to 0 $  and $v_a(M_a)\to \infty$ as $a \to \infty.$
\vskip .1 in

{\bf Proof:}  Suppose not.  Therefore, assume that $v_a$ is non-decreasing on $[0, R^{2-N}].$  

\vskip .1 in

\textbf{Claim:} For sufficiently large $a>0,$ there exists a $t_{a, \beta}\in (0, R^{2-N})$ such that $v_a(t_{a, \beta})=\beta$ and $t_{a, \beta}\to 0$ as $a\to \infty$, where $t_{a, \beta}$ is the smallest value of $t$ for which $v_a(t_{a, \beta})=\beta.$ 
\vskip .1 in
\textbf{Proof of Claim:} Rearranging terms and integrating $(\ref{vy000001})$ twice on $(0, t)\subset (0, J_{\frac{1}{\delta}, 1})\cap (0, R^{2-N}),$ we obtain
\begin{equation}
    v_a(t)+y(t) \int_0^t \frac{1}{y^2(t)} \int_0^s h(x)y(x)f(v_a(x)) \ dx \ ds= ay(t) \quad \forall \ t\in (0, J_{\frac{1}{\delta}, 1})\cap (0, R^{2-N}). \label{vy5}
\end{equation}
Now if $v_a<\beta$ for all $a,$ then $(\ref{vy5})$ implies
\begin{equation}
    v_a(t)\geq ay(t)=a A t^{\frac{1}{2}}J_{\frac{1}{\delta}}\left(\frac{2}{\delta}t^{\frac{\delta}{2}}\right)\quad \forall \ t\in (0, J_{\frac{1}{\delta}, 1})\cap (0, R^{2-N}). \label{vy9}
\end{equation}
Now let $0<t_0<J_{\frac{1}{\delta}, 1} \ \text{and} \ t_0<R^{2-N},$ then $ay(t_0)=a A t_0^{\frac{1}{2}}J_{\frac{1}{\delta}}\left(\frac{2}{\delta}t_0^{\frac{\delta}{2}}\right)\to \infty$ as $a\to \infty,$ so we see that $v_a$ exceeds $\beta,$ which contradicts our assumption. Therefore, there exists a smallest positive $t_{a, \beta}$ such that $v_a(t_{a, \beta})=\beta$ and $0<t_{a, \beta}<t_0.$ In addition, evaluating $(\ref{vy9})$ at $t_{a, \beta},$ we see
$$\frac{\beta}{a}\geq  A t_{a, \beta}^{\frac{1}{2}} J_{\frac{1}{\delta}}\left(\frac{2}{\delta}t_{a, \beta}^{\frac{\delta}{2}}\right).$$
It follows from this that $t_{a, \beta}\to 0$ as $a\to \infty.$ This proves the claim.
\vskip .1 in
Continuing now with our assumption that $v_a$ is non-decreasing on $[0, \ R^{2-N}],$ then we see that given $\epsilon>0$,  $\min\limits_{[\epsilon, \ R^{2-N}]} |v_a|=v_a(\epsilon)= \max\limits_{t\in [0, \ \epsilon]} v_a(t)\to \infty \  \text{as} \ a \to \infty$ by $(\ref{max1}).$ Thus $v_a\to \infty$ uniformly on $[\epsilon, \ R^{2-N}].$ Then for any $t\in [\epsilon, \ R^{2-N}]$, we obtain
\begin{align}
    \frac{h(t)f(v_a(t))}{v_a(t)}+\frac{1}{t^{2-\delta}}\geq \frac{h(t)f(v_a(t))}{v_a(t)} \geq h(t)v_a^{p-1}(t) \left ( 1-\left|\frac{g(v_a(t))}{v_a^p(t)}\right|\right)\notag \\ \geq h(R^{2-N})v_a^{p-1}(t) \left ( 1-\left|\frac{g(v_a(t))}{v_a^p(t)}\right|\right). \label{xyz01}
\end{align}
Since $v_a\to \infty$ uniformly as $a\to \infty \ \text{on}  \ [\epsilon, \ R^{2-N}],$ then $\left|\frac{g(v_a(t))}{v_a^p(t)}\right|\to 0$ uniformly as $a\to \infty$ by {\it (H1)}, and since $h$ is bounded from below on $[\epsilon, \ R^{2-N}],$ then the right-hand side of $(\ref{xyz01})$ goes to $\infty$ as $a\to \infty, \ \ie$
\begin{equation}
    \frac{h(t)f(v_a(t))}{v_a(t)}+\frac{1}{t^{2-\delta}}\to \infty \quad   \text{uniformly on} \ [\epsilon, \ R^{2-N}] \ \text{as} \ a\to \infty. \label{ift1}
\end{equation}
Now let $T_a:= \min\limits_{[\epsilon, \ R^{2-N}]} \left [ \frac{h(t)f(v_a(t))}{v_a(t)}+\frac{1}{t^{2-\delta}}\right].$ Then $T_a\to \infty \ \text{as} \ a\to \infty$ by $(\ref{ift1}).$  Next, let $\tilde y_a$ be the solution of
\begin{equation} \tilde y_a'' + T_a \tilde y_a=0  \label{y equation}  \end{equation}
with the same initial conditions as $v_a$ at $t= \epsilon$. Then $\tilde y_a$ is a linear combination of $\sin(\sqrt{T_{a}} t) $ and $\cos(\sqrt{T_{a}}t),$ and it is well known that $\tilde y_a$ has a zero on $\left[\epsilon, \epsilon+\frac{\pi}{\sqrt{T_a}}\right].$ Then for $a$ sufficiently large $\left[\epsilon, \epsilon+\frac{\pi}{\sqrt{T_a}}\right]\subset [\epsilon, \ R^{2-N}]$ and so $\tilde y_a$  has a zero, $z_{\tilde y_a}$,  on  $[\epsilon, \ R^{2-N}].$ Hence $\tilde y_a$ has a local maximum, $M_{y_a},$ on this interval with $\tilde y_a$ increasing on $[\epsilon, \ M_{\tilde y_a}]$.\\
Next, rewrite (\ref{e8}) as
\begin{equation} v_a'' + \left(\frac{h(t)f(v_a(t))}{v_a(t)}+\frac{1}{t^{2-\delta}}\right) v_a =0.  \label{muddy} \end{equation}
Multiplying   (\ref{muddy})  by $\tilde y_a,$ (\ref{y equation}) by $v_a$, and subtracting gives
$$  (v_a'\tilde y_a- v_a \tilde y_a')'   +  \left( \left(\frac{h(t)f(v_a(t))}{v_a(t)}+\frac{1}{t^{2-\delta}}\right) - T_a\right) v_a\tilde y_a = 0.$$
Integrating on $( \epsilon, M_{\tilde y_a})$ and recalling $\tilde y_a(M_{\tilde y_a})>0, $ $\tilde y_a'(M_{\tilde y_a})=0$ and $v_a'(M_{\tilde y_a})>0$ gives 
$$  \tilde y_a(M_{\tilde y_a})v_a'(M_{\tilde y_a}) + \int_{ \epsilon_m}^{M_{\tilde y_a}} \left( \left(\frac{h(t)f(v_a(t))}{v_a(t)}+\frac{1}{t^{2-\delta}}\right) - T_a\right) v_a \tilde y_a \, ds= 0.$$
However, this is impossible because the first  term on the left-hand side is positive, since  $\tilde y_a(M_{\tilde y_a})>0, $ $v_a'(M_{\tilde y_a})>0$,  and by the line after (\ref{ift1}) the second term on the left-hand side is also positive for $a $ sufficiently large. Thus we obtain a contradiction and so $v_a$ must have a  first local maximum, $M_a,$ on $[0, R^{2-N}]$ if $a$ is sufficiently large. It follows then from Lemma 2.5 that $v_a(M_a) \to \infty$ as $a \to \infty$. 

\vskip .1 in
Next, we show $M_a \to 0$ as $a \to\infty$.

\vskip .1 in
Suppose $M_a$ does not go to 0 as $a \to \infty$. Then, there is a $C_M>0$ such that $M_a > C_M$ for large $a$.
Again, since $t_{a,\beta}\to 0$ as $a\to \infty$ (by the claim at the beginning of this lemma), it follows that $t_{a,\beta}<\frac{M_{a}}{2}< M_a$ for large $a$ and from $(\ref{e7})$ we see that $h(t) f(v_a(t))+\frac{ v(t)}{t^{2-\delta}}>0,$ so $v_a''<0$ on $(t_{a,\beta}, M_a)$. Therefore,  using the concavity of $v_a$, we obtain
\begin{equation} v_a\left(\lambda_a t_{a, \beta} + (1-\lambda_a) M_a\right)  \geq  \lambda_a\beta + (1-\lambda_a)v_a(M_a) \quad  \text{for }0<\lambda_a<1.  \label{ringo} \end{equation}
Now choose $\lambda_a$ so that $\lambda_a t_{a, \beta} + (1-\lambda_a) M_a = \frac{M_a}{2}$ $\left( \ie  \ \lambda_a = \frac{M_a}{2(M_a-t_{a,\beta})}=  \frac{1}{2\left(1- \frac{t_{a,\beta}}{M_a}\right)}\right) $, then we see 
$\frac{t_{a,\beta}}{M_a}\leq \frac{t_{a,\beta}}{C_M} \to 0$ and so $\lambda_a \to \frac{1}{2}$  as $a \to \infty$. Then by (\ref{ringo}),
\begin{align} v_a\left(\frac{M_a}{2}\right) \geq  \lambda_a \beta + (1 - \lambda_a)v_a(M_a) =   \frac{\beta}{2(1-\frac{t_{a,\beta}}{M_a})}  + \left(1-\frac{1}{2(1-\frac{t_{a,\beta}}{M_a})}\right)v_a(M_a) \to \infty \ 
\textrm{ as } a \to \infty.  \label{tutti} 
\end{align}
Next, we  integrate (\ref{e7}) on $(t, M_a)$ with $\frac{M_a}{2} < t <M_a.$ Recalling that $v_a$ is increasing on $(\frac{M_a}{2}, M_a)$, and $v_a\to\infty$  on $(\frac{M_a}{2}, M_a)$ from (\ref{tutti}), it follows from  ({\it H1}) that $ f(v_a) \geq C_f v_a^p$ with $C_f>0$ on $(\frac{M_a}{2}, M_a)$, and we obtain
$$  v_a' = \int_{t}^{M_a}     h f(v_a)  \, ds \geq \int_{t}^{M_a} C_f h v_a^p   \, ds \geq  C_f v_a^{p}\int_{t}^{M_a} h  \, ds   \textrm{  \ \ on  } \left(\frac{M_a}{2}, M_a\right). $$
Thus, $$  \frac{ v_a'}{v_{a}^p} \geq C_f\int_{t}^{M_a} h \, ds \textrm{ \ \ on } \left(\frac{M_a}{2}, M_a\right).   $$
Integrating on $ (t,M_a)$ and recalling $p>1$, we obtain
\begin{equation}  v_{a}^{1-p}(t)  \geq v_a^{1-p}(t) - v_a^{1-p}(M_a) \geq (p-1)C_f \int_{t}^{M_a} \int_{s}^{M_a} h  \, dx \ ds  \quad \textrm{on  } \left(\frac{M_a}{2}, M_a\right). \label{1111} \end{equation}
Now, evaluate (\ref{1111}) at $t=\frac{M_a}{2}.$ Since $v_a\left(\frac{M_a}{2}\right)\to \infty$ as $a \to \infty$ from (\ref{tutti}) , then
$$  0 \leftarrow   v_a^{1-p}\left(\frac{M_a}{2}\right) \geq (p-1)C_f \int_{\frac{M_a}{2}}^{M_a} \int_{t}^{M_a} h \, ds \, dt \quad \text{as} \ a\to \infty.  $$
However,  as $a \to \infty$ the right-hand side goes to a positive  constant  (since $M_a \geq C_M>0$), and the left-hand side goes to 0, so we obtain a contradiction. Thus, $M_a \to 0$ as $a \to \infty$. This completes the proof.\qed

\vskip .1 in

{\bf Lemma 2.7:} Let $N>2$ and assume {\it (H1)}-{\it (H4)}. Suppose $v_a$ solves (\ref{e7})-(\ref{e8}). Then $v_a$ has a first zero, $z_a$, with $0< z_a <R^{2-N}$ if $a >0$ is sufficiently large. In addition, $z_a \to 0$ as $a \to \infty$.

\vskip .1 in

{\bf Proof:} From Lemma 2.6, we know that $v_a$ has a first local maximum, $M_a$ with $0<M_a<R^{2-N}$ for sufficiently large $a>0.$ Now suppose by way of contradiction that $v_a>0$ on $(M_a, R^{2-N})$ for all $a$ sufficiently large.\\
Then $\left (\frac{1}{2}  \frac{v_{a}'{^2}}{h} + F(v_a)\right)' = -\frac{v_a'{^2} h'}{2h^2}-\frac{ v_a v_a'}{h t^{2-\delta}}\geq 0 \  \textrm{on } (M_a, R^{2-N})$ when $v_a'\leq 0.$ It follows then that
\begin{equation} \frac{1}{2}  \frac{v_{a}'{^2}}{h} + F(v_a)  \geq F(v_a(M_a)) \textrm{ \ on } (M_a, R^{2-N}). \label{frutti} \end{equation}
$\textbf{Claim:}$ For sufficiently large $a>0,$ if $v_a>0$ on $(M_a, R^{2-N})$ then $v_a'<0$ on $(M_a, R^{2-N}).$
\vskip .1 in

$\textbf{Proof of Claim:}$ Suppose there is a first critical point $t_a$ with $M_a<t_a<R^{2-N},$ $v_a'(t_a)=0$ and $v_a(t_a)>0.$ Then from $(\ref{frutti}),$ we have $F(v_a(t_a))  \geq F(v_a(M_a))$. Since $v_a(M_a)\to \infty$ as $a\to \infty$ from Lemma 2.6, therefore $v_a(M_a)>\gamma$ for sufficiently large $a>0.$ So $F(v_a(M_a))>0$, and hence $F(v_a(t_a))>0$ for sufficiently large $a>0$ and thus $v_a(t_a)>\gamma.$ \\
On the other hand, since $v_a'<0$ on $(M_a, t_a),$ if $v_a''(t_a)< 0$ then $v_a$ has a local maximum at $t_a$ which contradicts that $v_a'<0$ on $(M_a, t_a).$ Thus $v_a''(t_a)\geq 0.$ Then from $(\ref{e7})$ we have $f(v_a(t_a))\leq 0.$ It follows then that $v_a(t_a)<\beta$, which contradicts that $v_a(t_a)>\gamma>\beta.$ Hence the claim is true.
\vskip .1 in
Next rewriting and integrating $(\ref{frutti})$ on $(M_a, R^{2-N})$ gives
\begin{equation}  \int_{M_a}^{R^{2-N}} \sqrt{h} \, ds \leq \int_{v_a(R^{2-N})}^{v_a(M_a)}  \frac{  dt}{ \sqrt{2} \sqrt{F(v_a(M_a)) - F(t)}}  
\leq \int_{0}^{v_a(M_a)}  \frac{  dt}{ \sqrt{2} \sqrt{F(v_a(M_a)) - F(t)}}. \label{diddley}  \end{equation}
Now, rewriting the right-hand side, letting $t= v_a(M_a)s$ and changing variables on the right-most integral in (\ref{diddley}), we obtain

$$ \int_{0}^{v_a(M_a)}  \frac{  dt}{ \sqrt{2} \sqrt{F(v_a(M_a)) - F(t)}}  = 
\frac{v_a(M_a)}{\sqrt{2}\sqrt{F(v_a(M_a))}} \int_{0}^1 \frac{ ds }{    \sqrt{1 - \frac{F(v_a(M_a)s)}{F(v_a(M_a))}       }         }    .  $$ 
We know from earlier that $v_a(M_a) \to \infty$ as $a \to \infty$ and so it follows from the superlinearity of $f$  (from {\it (H1)}) that  $\frac{v_a(M_a)}{\sqrt{2}\sqrt{F(v_a(M_a))}}   \to 0$,  and $\int_{0}^1 \frac{ ds }{    \sqrt{1 - \frac{F(v_a(M_a)s)}{F(v_a(M_a))}       }         }    
\to \int_{0}^{1} \frac{ ds }{\sqrt{1 - s^{p+1} }    } < \infty $ as $a \to \infty$.
Thus, it follows  that the right-hand side of (\ref{diddley}) goes to zero as $a \to \infty$ and therefore, we see
\begin{equation}
      \int_{M_a}^{R^{2-N}} \sqrt{h} \, ds \to 0  \textrm{ as }  a \to \infty \label{sqth1}.
\end{equation}

However, we know $M_a \to 0$ from Lemma 2.6, and so the left-hand side of  (\ref{diddley}) goes to  $ \int_{0}^{R^{2-N}} \sqrt{h} \, ds  >0$ as $a \to \infty$, which contradicts $(\ref{sqth1}).$ Thus, $v_a$ must have a first zero, $z_{a},$ with $M_a<z_a<R^{2-N}$ and $v_a'< 0$ on $(M_a, z_a].$\\
Now, rewriting  (\ref{frutti}) and integrating on $(M_a, z_{a})$ gives
\begin{equation}
    \int_{M_a}^{z_{a}} \sqrt{h}  \, ds  \leq \int_{0}^{v_a(M_a)}  \frac{  dt}{ \sqrt{2} \sqrt{F(v_a(M_a)) - F(t)}  } \label{hf001}.
\end{equation}
Again, since $f$ is superlinear, the right-hand side of $(\ref{hf001})$ goes to zero as $a\to \infty$, therefore $z_{a} - M_a \to 0$ as $a \to \infty.$ Since $M_a \to 0$ as $a \to 0$, by Lemma 2.6 it follows that $z_{a} = (z_{a} -M_a) + M_a \to 0$ as $ a \to \infty$.  \qed

\vskip .1 in

{\bf Lemma 2.8:} Let $N>2$ and assume {\it (H1)}-{\it(H4)}. Suppose $v_a$ solves (\ref{e7})-(\ref{e8}). Let $0< z_a <R^{2-N}.$ If $v_a(z_a)=0$, then $v_a'(z_a)\neq 0$ if $R >R_0$. In addition, $|v_a'(z_a)|\to \infty$ as $a\to \infty.$

\vskip .1 in

{\bf Proof:} Suppose $v_a(z_a)=0$ for some $0< z_a <R^{2-N}.$ Then there exists a local extremum on $(0, z_a)$, so suppose without loss of generality there is a local maximum, $M_a,$ with $M_a<z_a$ such that $v_a'(M_a)=0,$ and $v_a'(t)< 0$ for all $t\in (M_a, z_a].$ \\
From Lemma 2.3, we know that $|v_a(M_a)|> \gamma$ for $R>R_0$, so $F(v_a(M_a))>0.$ Since $\left (\frac{1}{2}  \frac{v_{a}'{^2}}{h} + F(v_a)\right)'\geq 0 \  \text{on} \ [M_a, z_a], \ \text{then} \ \frac{1}{2}  \frac{v_{a}'{^2}}{h} + F(v_a)$ is non-decreasing on $[M_a, z_a],$ and therefore
\begin{equation}
    0<F(v_a(M_a))\leq \frac{1}{2}  \frac{v_{a}'{^2}(z_a)}{h(z_a)}. \label{mnb01}
\end{equation}
 Hence from the above inequality, we see
$$|v_a'(z_a)|>0 \quad \text{for}\ R>R_0.$$
So, from Lemma 2.7 we know that $v_a(M_a)\to \infty$ and so from $(\ref{mnb01}),$ we see
$$ v_{a}'{^2}(z_a)\geq 2 h(z_a)F(v_a(M_a))\geq 2 h(R_0^{2-N})F(v_a(M_a))\to \infty \quad \text{as} \ a\to \infty.$$
Therefore, $|v_a'(z_a)|\to \infty$ as $a\to \infty.$ \qed

\vskip .1 in

We may now repeat this argument on $(z_a, R^{2-N})$ and show that $v_a$ has as many zeros as desired as $a \to \infty$ and at each zero, $z_a,$ we have $|v_a'(z_a)|\to \infty$ as $a\to \infty.$ 
\vskip .2 in

\section{Proof of Theorem 1}
First, assume $R>R_0$ where $R_0$ is obtained in Lemma 2.3.  Let $a>0.$ From Lemma 2.4,  we know that $v_a$ has at most a finite number of zeros on $(0, R^{2-N}).$ Now, let $n$ be a non-negative integer, and let
$$S_{n,R} = \{ a > 0  \ | \ v_a \textrm{ has exactly } n \textrm{ zeros  on } (0, R^{2-N}) \}.$$
From Lemma 2.4, we know some of the $S_{n,R}$ are non-empty so let $n_0$ be the smallest value of $n$ such that $S_{n_0,R}\neq \emptyset.$  It follows from Lemma 2.7 that $S_{n_0,R}$ is bounded from above. Let
$$a_{n_0}:=\sup S_{n_0,R}.$$
Now we prove that $v_{a_{n_0}}$ has $n_0$ zeros on $(0, R^{2-N}).$ For $a<a_{n_0},$ $v_a$ has exactly $n_0$ zeros on $(0, R^{2-N})$ and at each zero, $z,$ we have $v_a'(z)\neq 0$ by Lemma 2.8. Therefore, by continuity with respect to $a$ we see that $v_{a_{n_0}}$ has at least $n_0$ zeros on $(0, R^{2-N}).$ If $v_{a_{n_0}}$ has an $(n_0+1)^{st}$ zero on $(0, R^{2-N}),$ then so would $v_a$ for $a$ slightly smaller than $a_{n_0}$ (by continuity with respect to $a$) contradicting the definition of $a_{n_0}.$ Thus $v_{a_{n_0}}$ has exactly $n_0$ zeros on $(0, R^{2-N}).$
\vskip .1 in
Next suppose $v_{a_{n_0}}(R^{2-N})\neq 0.$ Without loss of generality suppose $v_{a_{n_0}}(R^{2-N})>0.$ Now for $a>a_{n_0}, \ v_a$ has an $(n_0+1)^{st}$ zero on $(0, R^{2-N})$ by definition of $a_{n_0}$ and so again by continuity with respect to $a$ it would follow that $v_{a_{n_0}}$ has an $(n_0+1)^{st}$ zero, a contradiction. Thus $v_{a_{n_0}}(R^{2-N})\leq0.$ Similarly, we can show that $v_{a_{n_0}}(R^{2-N})\geq0.$ Therefore $v_{a_{n_0}}(R^{2-N})=0$ and by Lemma 2.8, $v_{a_{n_0}}'(R^{2-N})\neq 0.$ 
\vskip .1 in
Now since $v_{a_{n_0}}'(R^{2-N})\neq 0,$ it follows that if $a$ is slightly larger that $a_{n_0}$ then $v_a$ has exactly $(n_0+1)$ zeros on $(0, R^{2-N}).$ Thus, $S_{n_0+1,R}$ is non-empty. Let
$$a_{n_0+1}:=\sup S_{n_0+1,R}.$$
In a similar way, we can show $v_{a_{n_0+1}}$ has exactly $(n_0+1)$ zeros on $(0, R^{2-N})$ and $v_{a_{n_0+1}}(R^{2-N})= 0, \ v_{a_{n_0+1}}'(R^{2-N})\neq 0.$   
\vskip .1  in
Continuing in this way, we see that for every $n\geq n_0$ there exists a solution of (\ref{e7})-(\ref{e8}) with exactly $n$ zeros on $(0, R^{2-N})$ and $v_{a_n}(R^{2-N})= 0.$ This completes the proof. \qed

\vskip .2 in

\section{Proof of Theorem 2}
Assume {\it (H1)-(H3)} and {\it (H5)} where $0< \alpha\leq 2.$ Suppose there is a $C^1[R, \infty)$ solution  $u$ of (\ref{DE1})-(\ref{DE100}). In particular, we assume $\frac{1}{|u|^q}$ is integrable. Set
$$E(r):=\frac{\frac{1}{2}u'{^2}(r)+\frac{ u^2(r)}{2r^{2+(N-2)\delta}}}{K(r)}+F(u(r)).$$
\textbf{Claim 1:} Let $u$ be a $C^1[R, \infty)$ solution of (\ref{DE1})-(\ref{DE100}). Then $E'\leq 0$ and $E(r)\to 0$ as $r\to \infty.$

\vskip .1  in

\textbf{Proof of Claim 1:} A straightforward calculation shows
\begin{equation}
    E'(r)=-\left [2(N-1)+\frac{r K'(r)}{K(r)}\right]\frac{u'{^2}(r)}{2rK(r)}-\left [2+(N-2)\delta+\frac{rK'(r)}{K(r)}\right ]\frac{ u^2(r)}{2r^{3+(N-2)\delta}}. \label{ed}
\end{equation}
By {\it (H5)} we have $\alpha+\frac{rK'}{K}> 0$ so since $\delta>0$ and $0<\alpha\leq 2$ it follows that $2+(N-2)\delta+\frac{rK'}{K}\geq \alpha+\frac{rK'}{K}> 0$ and $2(N-1)+\frac{r K'}{K}\geq \alpha+\frac{rK'}{K}> 0,$ therefore $E'(r)\leq0$ on $[R, \infty),$ so $E$  is non-increasing on $[R, \infty)$ when $0<\alpha\leq 2.$ Also $E\geq -F_0$ from $(\ref{boss2}),$  so $E$ is bounded from below. Therefore, there exists $L$ such that 
$$ \lim\limits_{r\to \infty} E(r)=L.$$
In addition, since $u\to 0$ as $r\to \infty$ implies $F(u)\to 0$ and also $\frac{\frac{ u(r)^2}{r^{2+(N-2)\delta}}}{K(r)}\leq \frac{ u(r)^2}{r^{(2-\alpha)+(N-2)\delta}} \to 0 \ \text{as} \ r\to \infty \ (\because \ 0<\alpha\leq 2 \ \text{and} \ (\it{H5})),$ and since $\frac{u'{^2}}{2K}\geq 0,$ therefore $\lim\limits_{r\to \infty} \frac{u'{^2}}{2K}=L,$ and hence $L\geq 0.$ 
\vskip .1 in
Now, we show that $L=0.$ Suppose not. Suppose $L>0,$ so then $\lim\limits_{r\to \infty} \frac{u'{^2}}{2K}=L>0,$ and thus there exists $r_L>0$ such that
\begin{equation}
    \frac{u'{^2}}{2K}>\frac{1}{2}L>0 \quad \forall \ r\geq r_L. \label{non01}
\end{equation}
So $u'(r)\neq0$ for all $r\geq r_L.$ Without loss of generality, assume $u'(r)<0$ for all $r\geq r_L.$ Now taking the square root in $(\ref{non01})$ and using $(\it{H5}),$ we see
$$-u'(r)>\sqrt{\frac{L}{2}}C_{10} r^{-\frac{\alpha}{2}}\quad \forall \ r>r_L, \ \text{and for some} \ C_{10}>0.$$
Integrating on $(r_L, r),$ we obtain
$$0\leftarrow u(r)\leq \frac{\sqrt{\frac{L}{2}}C_{10} r_L^{1-\frac{\alpha}{2}}}{1-\frac{\alpha}{2}}+u(r_L)-\frac{\sqrt{\frac{L}{2}}C_{10} }{1-\frac{\alpha}{2}}r^{1-\frac{\alpha}{2}}\to -\infty \quad \text{as} \ r\to \infty,$$
which is a contradiction. This proves the claim.

\vskip .1  in

\textbf{Claim 2:} Let $u$ be a $C^1[R, \infty)$ solution of (\ref{DE1})-(\ref{DE100}). If $z\in [R, \infty)$ and $u(z)=0$, then $u'(z)\neq0.$

\vskip .1  in

\textbf{Proof of Claim 2:} Suppose by way of contradiction that $u(z)=u'(z)=0.$ Then $E(z)=0.$ From Claim 1, we know that $E$ is non-increasing and $E(r)\to 0$ as $r\to \infty,$ therefore $E(r)= 0$ for all $r\geq z.$ Hence $E'(r)= 0 \ \text{for all} \ r\geq z$, and this with (\ref{ed}) gives
$$0\geq -\left[2(N-1)+\frac{r K'(r)}{K(r)}\right]\frac{u'{^2}(r)}{2rK(r)}=\left[2+(N-2)\delta+\frac{rK'(r)}{K(r)}\right]\frac{ u^2(r)}{2r^{3+(N-2)\delta}}\geq 0,$$
which forces
\begin{equation}
  u(r)=u'(r)=0 \quad \forall \ r\geq z. \label{rz} 
\end{equation}
But then $\frac{1}{|u|^q}$ is not integrable on $[R, \infty)$ and so $u$ is not a solution. Thus the claim is true.

\vskip .1  in

\textbf{Claim 3:} Let $u$ be a $C^1[R, \infty)$ solution of (\ref{DE1})-(\ref{DE100}). Then $u$ can have at most a finite number of zeros.
\vskip .1  in

\textbf{Proof of Claim 3:} Suppose not. Let $u$ have infinitely many zeros, say $z_k, \ \text{with} \ z_k<z_{k+1}  \ \forall \ k\in \mathbb N$. Then $u$ must have infinitely many extrema, say $m_k,$ such that 
$$z_k<m_k<z_{k+1}, \ \forall \ k\in \mathbb N.$$
{\bf Case 1:} $z_k\to z^*$  for some $R<z^*<\infty$
\vskip .1  in
Then $u(z^*)=0$ and $u'(z^*)=0$ $(\because \ m_k\to z^* \ \text{and } \ u(m_k)=0 \ \forall \ k),$ which contradicts Claim 2.

\vskip .1  in
{\bf Case 2:} $z_k\to \infty$ 
\vskip .1  in
In this case, $m_k\to \infty.$ Then by Lemma 2.3, we obtain
$$u(m_k)\geq C_7>0,$$
but this contradicts that $u(r)\to 0 \ \text{as} \ r\to \infty.$ This proves the claim.

\vskip .1  in
We continue now with the proof of Theorem 2.
\vskip .1  in
We know that $u$ can have at most finitely many zeros. So, suppose $u$ has $n$ zeros $z_i \ (i=1, 2,\cdots, n)$ such that $z_i<z_{i+1}, \ i=1, 2, 3,\cdots, n-1.$ Since $u'(z_n)\neq 0$ by Claim 2, therefore without loss of generality, we may assume that $u(r)\geq 0 \ \text{for all } \ r\geq z_n.$ Then $u$ must be monotone for large $r$. Otherwise, $u$ would have infinitely many extrema, $m_k$ with $m_k\to \infty$ and by Lemma 2.3 we obtain 
$$u(m_k)\geq C_7>0 \quad \text{for positive local maxima} \ m_k.$$
Since $m_k\to \infty,$ then this contradicts that $u(r)\to 0 \ \text{as} \ r\to \infty.$ Thus $u$ is monotone for large $r.$ 
\vskip .1  in
Now, since $u(r)\to 0 \ \text{as} \ r\to \infty $ and $u$ is monotone for large $r$, therefore without loss of generality, we may assume that there is an $r_0>R$ such
\begin{equation}
    0\leq u(r)\leq \frac{\beta}{2} \ \text{and} \ u'(r)\leq 0 \quad \forall \ r>r_0.\label{du10}
\end{equation}
Next, we multiply (\ref{DE1}) by $r^{1+\frac{\delta}{2}}u^q, $ integrate by parts on $(r_0,r)$   and  obtain
\begin{equation} 
r^{1+\frac{\delta}{2}}u^q u' \geq C_{11}   - \int_{r_0}^r  t^{1+\frac{\delta}{2}}K(t) u^q(t) f(u(t)) \ dt,  \label{wrr} 
\end{equation}
where $C_{11}:=r_0^{1+\frac{\delta}{2}}u^q(r_0) u'(r_0)  + \frac{u^{q+1}(r_0) r_0^{\frac{\delta}{2}-(N-2)\delta}}{\frac{\delta}{2}-(N-2)\delta}.$ For large $r$ we have $u\to 0$, so it follows from {\it (H2) } and (\ref{du10}) that  $- u^q f(u) > \frac{1}{2}$ for large $r$. Thus, from (\ref{wrr}), we obtain
\begin{equation}
u^q u' \geq  \frac{C_{11}}{r^{1+\frac{\delta}{2}}}  +  \frac{1}{2r^{1+\frac{\delta}{2}}}  \int_{r_0}^r  t^{1+\frac{\delta}{2}}K(t)  \, dt.  \label{bach}  
\end{equation}
Integrating again on $(r_0, r)$ gives
\begin{equation}
u^{q+1} \geq u^{q+1}(r_0) + \frac{2C_{11}(q+1)}{\delta}(r^{-\frac{\delta}{2}}+r_0^{-\frac{\delta}{2}})+ \frac{q+1}{2} \int_{r_0}^r \frac{1}{t^{1+\frac{\delta}{2}}}\int_{r_0}^t s^{1+\frac{\delta}{2}} K(s) \, ds \, dt. \label{daddo}  
\end{equation}
Now, recalling {\it (H5)} we see that $r_0^\alpha K(r_0)\leq r^\alpha K(r)$ for $r\geq r_0.$ Then we estimate the integral in (\ref{daddo}) and we obtain
\begin{align}
\int_{r_0}^r \frac{1}{t^{1+\frac{\delta}{2}}}\int_{r_0}^t s^{1+\frac{\delta}{2}} K(s) \, ds \, dt \geq  r_0^\alpha K(r_0)  \int_{r_0}^r \frac{1}{t^{1+\frac{\delta}{2}}}\int_{r_0}^t s^{1+\frac{\delta}{2}-\alpha} \, ds \, dt  \notag\\ =\frac{r_0^\alpha K(r_0)}{2+\frac{\delta}{2}-\alpha}\int_{r_0}^r \frac{1}{t^{1+\frac{\delta}{2}}} (t^{2+\frac{\delta}{2}-\alpha} - r_0^{2+\frac{\delta}{2}-\alpha}) \, dt. \label{aunt ruth} 
\end{align}
If $0< \alpha <2$, then it follows from (\ref{aunt ruth}) that
\begin{align}
\int_{r_0}^r \frac{1}{t^{1+\frac{\delta}{2}}}\int_{r_0}^t s^{1+\frac{\delta}{2}} K(s) \, ds \, dt  \geq \frac{r_0^\alpha K(r_0)}{2+\frac{\delta}{2}-\alpha}\left( \frac{ r^{2-\alpha} - r_0^{2-\alpha}}{2-\alpha} + \frac{2r_0^{2+\frac{\delta}{2}-\alpha}}{\delta}(r^{-\frac{\delta}{2}}-r_0^{-\frac{\delta}{2}})\right) \to \infty \textrm{ as }  r \to \infty. \label{mozart}\end{align}
If $\alpha =2$, then it follows from (\ref{aunt ruth})  that
\begin{align}
\int_{r_0}^r \frac{1}{t^{1+\frac{\delta}{2}}}\int_{r_0}^t s^{1+\frac{\delta}{2}} K(s) \, ds \, dt \geq   \frac{r_0^\alpha K(r_0)}{2+\frac{\delta}{2}-\alpha} \left(\ln(r) - \ln(r_0) - \frac{2r_0^{\frac{\delta}{2}}}{\delta}(r^{-\frac{\delta}{2}}-r_0^{-\frac{\delta}{2}})\right) \to \infty \textrm{ as }  r \to \infty.  \label{beethoven} \end{align}
Now, notice that the first two terms on the right-hand side of (\ref{daddo}) are bounded,  so it follows from (\ref{mozart})-(\ref{beethoven}) that the right-hand side of (\ref{daddo}) goes to infinity as $r\to \infty$, but by assumption the left-hand side of (\ref{daddo}) goes to zero as $r \to \infty$. Thus, we obtain a contradiction, so the theorem is proved. \qed

\end{document}